\documentclass[12pt]{amsart}
\usepackage{amsmath,amsfonts,amssymb,amsthm}
\usepackage[mathscr]{eucal}
\usepackage{enumitem,kantlipsum}
\usepackage{color}

\usepackage{datetime2}
\usepackage[numbers,sort&compress]{natbib}
\usepackage[dvipsnames]{xcolor}
\usepackage{graphicx}
\usepackage[autostyle]{csquotes}
\usepackage[abs]{overpic}
\usepackage{caption}
\RequirePackage[colorlinks=true, linkcolor = red,citecolor=blue,urlcolor=blue,backref=page]{hyperref}
\usepackage{float}

\usepackage[normalem]{ulem}
\usepackage[numbers]{natbib}
\numberwithin{equation}{section}
\theoremstyle{definition}
\newtheorem{definition}{Definition}[section]
\theoremstyle{definition}

\theoremstyle{plain}
\newtheorem{theorem}[definition]{Theorem}

\newcommand{\nc}{\newcommand}
\nc{\I}{{\mathbf 1}}
\newcommand{\remove}[1]{}
\nc{\bN}{{\mathbf N}}
\nc{\cB}{{\mathcal B}}
\nc{\cF}{{\mathcal F}}
\nc{\cX}{{\mathcal X}}
\nc{\cS}{{\mathcal S}}
\nc{\cC}{{\mathcal C}}
\nc{\cP}{{\mathcal P}}
\nc{\cN}{{\mathcal N}}
\nc{\R}{{\mathbb R}}

\nc{\E}{{\mathbb E}}
\nc{\C}{{\mathbb C}}

\nc{\N}{{\mathbb N}}
\nc{\Z}{{\mathbb Z}}
\nc{\BX}{{\mathbb X}}
\nc{\BM}{{\mathbb M}}
\nc{\BY}{{\mathbb Y}}
\nc{\bfW}{{\mathbf W}}
\nc{\bfY}{{\mathbf Y}}
\nc{\bfZ}{{\mathbf Z}}
\nc{\bx}{\mathbf{x}}
\nc{\bd}{\partial}
\def\1{\mathbf{1}}

\nc{\BP}{\mathbb{P}}
\nc{\BE}{\mathbb{E}}
\nc{\BQ}{\mathbb{Q}}

\newcommand{\beas}{\begin{eqnarray*}}
\newcommand{\eeas}{\end{eqnarray*}}
\newcommand{\bes} {\begin{equation*}}
\newcommand{\ees} {\end{equation*}}
\newcommand{\be} {\begin{equation}}
\newcommand{\ee} {\end{equation}}
\newcommand{\bea} {\begin{eqnarray}}
\newcommand{\eea} {\end{eqnarray}}

\emergencystretch=\maxdimen
\makeatletter
\g@addto@macro\bfseries{\boldmath}
\makeatother

\begin{document}

\title{The Life and Works of K.G Ramanathan}
\author{Maitreyo Bhattacharjee }
\address{School of Mathematical and Computational Science,  Indian Association for the Cultivation of Science,  Kolkata.}
\email{ug2020mb2353@iacs.res.in,  maitreyomaths@gmail.com}

\subjclass[2020]{ \textbf{Primary} 01A70, 
11A07 
11A25 
\textbf{Secondary} 01A75, 
01A73;  	
01A72. 
}
\date{October 25, 2022}

\maketitle

\begin{abstract}
    Prof. K.G Ramanathan was a legendary Indian Mathematician, working in Number Theory and a prolific Institution builder. Apart from this, he was an excellent teacher and influenced several brilliant students. In this article, we overview his life and discuss some of his significant mathematical contributions.
\end{abstract}

\section{Early life and education}

\begin{figure}[h]
\centering
\includegraphics[width=6 cm,height=5 cm]{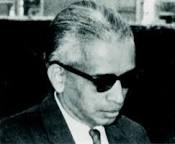}
\medskip
\caption{K. G Ramanathan (1920-1992)}
\end{figure}

Kollagunta Gopalaiyer Ramanathan, or popularly known to his colleagues,
friends and students as KGR, was a distinguished Indian mathematician, well known for his outstanding work in number theory and  contributions in building a strong school of mathematics research and teaching in India. His influence on the mathematical scene in India post-independence was huge. Ramanathan was born on the 13th of November, 1920 in Hyderabad, a city which is situated in southern India, and currently the capital of Telangana. His mother, Ananthalaxmi, died at an early age. His father’s name was Kollagunta Gopal Iyer. Apart from him, there were two sisters and one brother in his family. Ramanathan’s schooling was done from Wesleyan Mission High School, Secunderabad. He studied at the Nizam College, Hyderabad \cite{NEW} and obtained his undergraduate degree from Osmania University in 1940 and his Master’s Degree from the University of Madras in 1942, studying at the Loyola College. He also worked as a Research Scholar and taught for a while in the latter before moving to the United States for higher studies.  During his stay at Madras, Ramanathan came in contact with many famous mathematicians like Professors Vaidyanathaswamy and Vijayaraghavan, who inspired him. His first research paper was on Demlo Numbers \cite{KDN}, and it was published in 1941. During this period, he published some other papers as well, in reputed journals like the Mathematics Student and Journal of the Indian Mathematical Society.
\par In 1951, he received his Doctoral Degree in Mathematics from Princeton University, where his advisor was none other than \textbf{Emil Artin}, considered as one of the co founders of the theory of modern abstract algebra (along with Emmy Noether). The title of his doctoral thesis was “The Theory of Units of Quadratic and Hermitian Forms.” At Princeton, he came in contact with legendary mathematicians like Hermann Weyl and Carl Siegel, whose mathematics heavily influenced his own research later. When he was the Institute for Advanced Study, Ramanathan assisted Prof. Weyl, one of the most significant Mathematicians and Theoretical Physicists of the last century, who was a Member of the Institute at that time. The mathematics he learnt during these crucial early years had a lasting effect on his entire career.

\section{Returning back to India}

\par After completing
his PhD in the USA, Ramanathan returned to India in 1951 and joined K. Chandrasekharan at the
School of Mathematics, Tata Institute of Fundamental Research (TIFR), Mumbai. The
department had just started its journey. The decade of 1950s is very crucial for Indian
mathematics, as it was during that time when Ramanathan and Chandrasekharan, with
full support and encouragement from Homi Jehangir Bhabha, FRS (the first Director of the Institute),
were working on transferring the department into a notable centre for mathematics.
They attracted some amazing research scholars and also invited many leading
mathematicians like C.L Siegal, Samuel Eilenberg, Laurent Schwartz and Kiyosi Ito to visit the Department and deliver courses in many areas of advanced mathematics. For example, Eilenberg gave a course on the recent developments in Algebraic Topology in 1955, and Schwartz lectured on complex analytic manifolds. \cite{NNL} The work culture in those days was such that the students were given complete academic freedom to read and explore whatever they wanted, as long as Ramanathan and Chandrasekharan themselves were convinced that the student was serious enough. They however, asked difficult questions as Interview panel members during the time of admission of PhD students. With
his unending enthusiasm for good mathematics and dedication towards the
development of teaching and research of math in our country, Ramanathan was
successful in setting up a strong Number theory school at TIFR. He had amazing international connections, by attending various International Conferences and Meetings. For example, he attended the landmark and immensely successful \textbf{Tokyo Nikko Conference} on Algebraic Number Theory in September 1955 (where his advisor Artin was also present), and gave a talk titled \enquote{Units of fixed points in involutorial algebras}\cite{PTK}. The event was the first major Mathematics event held in Japan after the Second World War, and is special for many reasons. It attracted 77 participants, which included legendary number theorists from the West, like Artin, Chevalley, Deuring, Serre, Weil, as well as the new generation of Japanese mathematicians including Iwasawa, Satake, Shimura, Taniyama, Takashi Ono (father of the famous mathematician Ken Ono). At an after dinner impromptu lecture(which was not a part of the official program), Weil told the inspiring story of Ramanujan to the Japanese mathematicinas, who were not familiar with him \cite{WR}. It was at this Conference that the Japanese mathematician Yutaka Taniyama first stated the famous Taniyama-Shimura-Weil Conjecture, about elliptic curves and modular forms. \cite{KFT}
After almost 3 decades, Ken Ribet \cite{EC} made a breakthrough by proving the Epsilon Conjecture, first proposed by Serre, thereby establishing that the Taniyama-Shimura-Weil Conjecture implies the Fermat's Last Theorem, and the history. The work at the 1955 Conference, has an important role to play in Wiles's famous proof. Thus, the Tokyo 1955 Conference has an important place in the history of Number Theory. 

\begin{figure}[h]
\centering
\includegraphics[width=11.5 cm,height=8 cm]{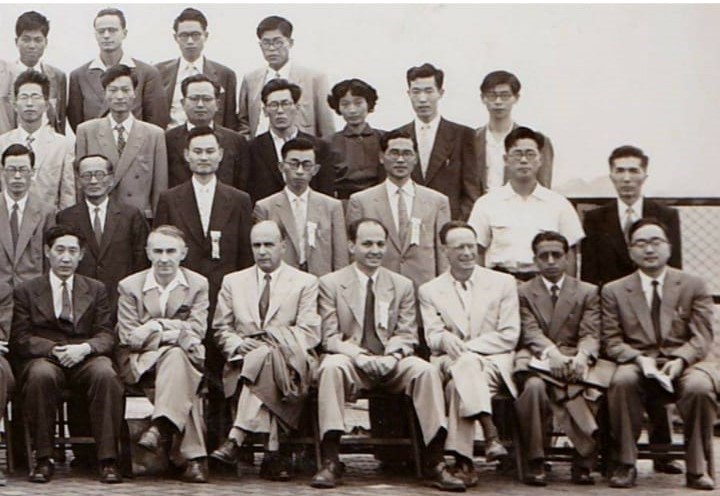}
\medskip
\caption{\small{Part of the group photo at the 1955 Tokyo Nikko Conference. \textbf{Top} Row, from left - Morikawa, Serre, Taniyama Yamazaki; \textbf{Third} Row, from left - Tsuzuku, Takashi Ken Ono, Azumaya, Nagai, Keiko Satake, Ichiko Satake, Hayashida; \textbf{Second} Row, from left - Shimura, Shinziro Mori, Hitotumatu, Terada, Nakano, Kawai, Ishii; \textbf{Seated}, from left - Iyanaga, Chevalley, Brauer, Zelinsky, Weil, KGR, Iwasawa. (taken from Photo Key in \cite{PTK}) \\ Courtesy - \textbf{My Parents' Generation}, Ken Ono and Amir D. Aczel \cite{KPGB} }.}
\end{figure}

Today, the School of
Mathematics, Tata Institute of Fundamental Research (TIFR), Mumbai is a leading centre
of pure mathematics research with brilliant faculty and international reputation, and this
would never have been possible without the effort of the two. A visionary, Ramanathan
realized the need for establishing a centre for Applied Mathematics in India. The idea of
establishment of a joint TIFR-IISc programme (which is the TIFR Centre for Applicable
Mathematics (CAM)) to be operated from the campus of the Indian Institute of Science,
Bangalore in 1975 was his brainchild. In fact, he preferred the phrase \enquote{applicable} to \enquote{applied}. Today, the TIFR-CAM has become an important Institution, where faculties are actively working on many aspects of differential equations and numerical
methods. In February 2005, an Indo-French Workshop on Partial Differential Equations and their Applications was organized at the Department of Mathematics, IISc Bangalore \cite{WEB}. Several distinguished speakers form both the Indian as well as French math community participated in it. This event was dedicated to the memory of KGR and \textbf{Prof. Jacques-Louis Lions} (father of Pierre Louis Lions, Fields Medal, 1994), who were pioneers in initiating the Indo-French cooperation in Applied Mathematics in the 1970s. Thus, Ramanathan holds the distinction of leading the foundation of both pure
as well as applied mathematics schools in India \cite{KBA}. Such excellent initiative can indeed be compared with that of Prasanta Chandra Mahalonbis and C.R Rao in setting up the renowned school of Statistics and Probability at  the Indian Statistical Institute, Kolkata and later in its branches in Bengaluru and Delhi. Almost 7 decades have passed after the
establishment of the TIFR, and today we have so many other good places to teach and do mathematical research in India. During his time as at the Tata Institute in Mumbai, Ramanathan closely
and enthusiastically interacted with numerous research scholars during their long walks
along the Arabian Sea, and with his vast knowledge, exposed them to many deep areas
of pure mathematics, which were still not that much popular in India. He also had a nice sense of humour and had the
habit of sharing many interesting anecdotes about other mathematicians among his
students. An excellent teacher and expositor of the subject, Ramanathan’s lessons and guidance had a
great impact on their future mathematical career. Some of the students went on to become acclaimed
mathematicians, like MS Narasimhan, C S Seshadri and R Sridharan \cite{RSK} (all three of them were Shanti Swarup Bhatnagar Awardees). An interesting anecdote is worth making in this regard. After completing their Doctoral degree in India (where they were mostly self taught), both Narasimhan and Seshadri had went for Postdoctoral stints in Paris, with the former being mentored by Laurent Schwartz (Fields Medal, 1950) and the latter by Claude Chevalley \cite{NN}. Even before they went to Paris, it was KGR who drew their attention to a very important paper of Weil on vector bundles on compact Riemann surface (which did not directly belong to Ramanathan's own area of research), titled “Generalisation des fonctions abeliennes". KGR was himself made aware of this work of Weil by Seigel. This was played a very crucial role in their research, as they got influenced and began working on the topic. In 1964, they finally established their celebrated result, the \textbf{Narasimhan Seshadri Theorem}, which concerns the stability of holomorphic vector bundles over a Riemann Surface.

\begin{figure}[h]
\centering
\includegraphics[width=9.5 cm,height=6.5 cm]{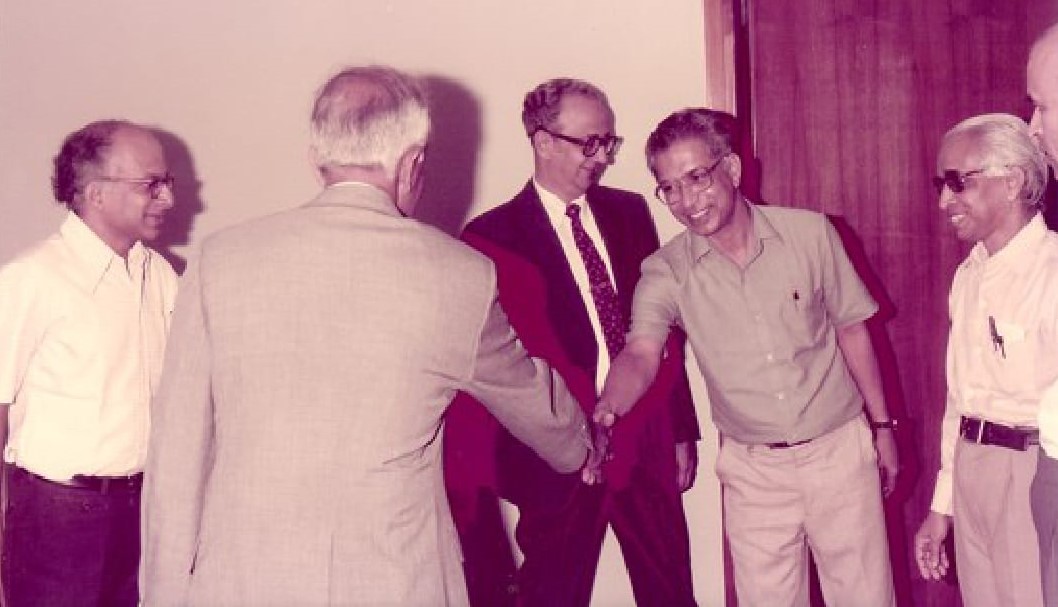}
\medskip
\caption{\small{KGR (second from right) with (from left) M.S Raghunathan, S. Chandrasekhar(Nobel Laureate in Physics, 1983), B.V Sreekantan, M.S Narasimhan and K. Ramachandra, during Chandrasekhar's 1987 visit to the Institute. \\ Courtesy - \textbf{TIFR Archive}}.}
\end{figure}

\section{Mathematical Works}
Now we would discuss some mathematical works of Ramanathan. Broadly speaking, his interests were primarily in Number Theory and related fields. He was extremely interested in exploring and extending the published and unpublished works of
Ramanujan and had a few publications on the same. He was also actively involved in
studying his Notebooks, and wrote some papers \cite{CPR},   after being motivated by the
beautiful mathematics of Ramanujan, which were mainly based on the topics of congruence properties of
some arithmetical functions, Ramanujan’s Trigonometric Sums and Ramanujan type
identities. One of his remarkable papers was “Identities and
Congruence of the Ramanujan Type” \cite{PPP}, published in 1950. We
know that Ramanujan was famous for his work on $P(n)$, which denotes the number of
“unrestricted \textbf{partitions}” of a natural number n. For example, $P(5)=7$, i.e., there are 7
partitions of 5, namely - (4+1), (3+2), (3+1+1+1), (2+2+1), (2+1+1+1) and (1+1+1+1+1).
One may notice that we are trying to write $n$ as sum of positive integers $\leq~n$, where the
summands are in \textbf{non-increasing} order. Partition is an important object of study in Number Theory. Below, we state one of Ramanathan’s result from that paper: \\
\begin{theorem}
Let $v>0$,
\[ \sum_{n=0}^{\infty} P_v(n)x^n = \frac{1}{[(1-x)(1-x^2) \cdots]^{-v}} \]
The coefficients $P_v(n)$  satisfy the following properties :
\begin{enumerate}
    \item If 24m $\equiv v (\text{mod}~ 5^a)$, then ~ $P_v(m) \equiv 0$ :
    \begin{center}
    \begin{enumerate}
        \item $(mod~5)$ if v $\equiv 12,17,22,27~ (mod~30)$ \\
         \item ($mod~5^{\big[ \frac{a-1}{2} \big]}$) if $v \equiv 15,20,25~ (mod~30)$ \\
        \item ($mod~5^{\big[ \frac{a}{2} \big]}$) if $v \equiv 3,4,8,9~ (mod~30)$ \\
        \item ($mod~5^{\big[ \frac{a+1}{2} \big]}$) if $v \equiv 16,21,26~ (mod~30)$ \\
        \item ($mod~5^{a-1}$) if $v \equiv 0,5,10~ (mod~30)$ \\
        \item ($mod~5^{\big[ \frac{a+2}{2} \big]}$) if $v \equiv 2,7~ (mod~30)$ \\
        \item ($mod~5^{a}$) if $v \equiv 1,6,11~ (mod~30)$
        
    \end{enumerate}
    \end{center}
   \item If 24m $\equiv v (\text{mod}~ 7^a)$, then ~ $P_v(m) \equiv 0$ :
   \begin{center}
    \begin{enumerate}
        \item $(mod~7)$ if v $\equiv 8,15,22~ (mod~28)$ \\
         \item ($mod~7^{\big[ \frac{a-1}{2} \big]}$) if $v \equiv 14,21~ (mod~28)$ \\
        \item ($mod~7^{\big[ \frac{a}{2} \big]}$) if $v \equiv 2,3,5,6~ (mod~28)$ \\
        \item ($mod~7^{\big[ \frac{a+1}{2} \big]}$) if $v \equiv 11,18,25~ (mod~28)$ \\
        \item ($mod~7^{a-1}$) if $v \equiv 0,7~ (mod~28)$ \\
        \item ($mod~7^{\big[ \frac{a+2}{2} \big]}$) if $v \equiv 1~ (mod~28)$ \\
        \item ($mod~7^{a}$) if $v \equiv 4~ (mod~28)$
        
    \end{enumerate}
    \end{center}
\end{enumerate}

\end{theorem}

Here, [.] denotes the usual “greatest integer function.” 
The proof of Ramanathan uses some techniques developed earlier by mathematicians H.
Rademacher (specially construction of some special functions) and G.N Watson, and
uses the Dedekind Modular form, defined as follows :
\[ \tau(\eta)=e^{\frac{i\pi\tau}{12}} \prod_{n=1}^{\infty} (1-e^{2n\pi i \tau}) \]
and the “hauptmodul” $\phi(\tau)$ of a modular function  The proof also rests on \textbf{8
lemmas}, where some new entire modular functions are defined and manipulated
accordingly. In particular, the following function is constructed :
\[ F_{n,\nu}=\phi(\tau)^{kt}\sum_{i>0}a_i(n)p^{i-1}\phi(\tau)^i \]
(the $a_i(n)$ being integers depends on ,$n,p$ and $\nu$) used with other functions for proving
the congruence properties for $P_{\nu}(n)$. \\ 
\par Ramanujan had conjectured that :
\[ P(n) \equiv 0 ~(\text{mod}~ 5^a7^b11^c)\]
It is not known whether this is true for all c. Though he did not prove these, he indicated that these congruences may be deduced from identities of this kind :
\[ P(4)+P(9)x+\cdots=5 \cdot \frac{[(1-x^5)(1-x^{10})\cdots]^5}{[(1-x)(1-x^2)\cdots]^6}\]
The above theorems established by Ramanathan for are analogous to these
identities of Ramanujan, and can be applied to study them in details, and perhaps derive
some similar formulas.  Ramanujan considered the following sum in one of his papers published in 1918 \cite{RTS} :
$$ c_q(n) = \sum_{\substack{1 \leq a \leq q \\ (a,q)=1}} e^{\frac{2\pi ian}{q}} $$
$c_q(n)$, usually called the Ramanujan Sum, is a function of two positive integer variables q and n. Ramanathan studied this
function in great details \cite{RRR1, RRR2} and found many new properties. Now, if a natural number n has the following property :
\[ n \equiv e_1 + e_2 + \cdots (\text{mod}~ m)\]
then n is said to be relatively partitioned (mod m). von Sterneck used a certain two
variable function $f(m,n)$ to study and write many formulas related to this theory in
closed form. Ramanathan showed that $f(m,n)$ is same as that of the Ramanujan Sum $c_q(n)$. He also found many other applications of Ramanujan’s Trigonometric sum in the theory developed by von Sterneck. He had a famous paper \cite{KLF}  on
the applications of Kronecker’s limit formula. Ramanathan also did important work in
other areas of Number Theory, namely Diophantine Inequalities and automorphic
functions related to Siegel’s formula. He also generalised certain results of
mathematicians Hecke and Maass. KGR also extensively worked on the
properties of unit groups of quadratic and Hermitian forms over algebraic number fields \cite{KAJ, KAJ1} with one of his papers coming out in the $\textit{Annals of Mathematics}$.
 In modern number theory, arithmetic groups are an important class of objects, which are of interests to Topologists as well. It was KGR who popularized the study of Arithmetic Groups in India, and later it became one of the main research topics at the Tata Institute. His work on the analytic and arithmetic theory of quadratic forms over involutorial division algebra is internationally recognized. KGR, along with Raghavan also worked and contributed significantly towards the famous Oppenheim Conjecture \cite{RR2}, which concerns quadratic forms in many variables. The conjecture was finally settled in the affirmative by Margulis in 1987. Coauthoring a paper with the Indo-Canadian mathematician Mathukumalli V. Subbarao \cite{ERK} gave him an Erdős number of \textbf{2}.
 
 Interested readers may check out an obituary written by Prof. S. Raghavan \cite{SRO} where some other aspects of his works (mainly involving discrete groups) are covered. It also contains the complete list of \textbf{48} research papers and \textbf{2} books of Ramanathan.
\section{Awards, Honours and Recognition}
\par Throughout his career, Ramanathan received numerous honours and recognition. He was
awarded the coveted Shanti Swarup Bhatnagar Prize in Mathematical Science (one of
the highest recognitions for an Indian researcher) in the year 1965, the Padma Bhushan,
the third highest civilian award awarded by the Government of India in 1983 and the
Homi Bhabha Medal of the Indian National Science Academy (INSA) in 1984.
Ramanathan was the Fellow of the Indian Academy of Sciences, Bangalore, Indian
National Science Academy, New Delhi, the Jawaharlal Nehru Fellowship and the
Foundational Fellow of the Maharashtra Academy of Sciences. He was on the Editorial
Board of the prestigious journal Acta Arithmetica for almost three decades, and was the Editor of the Journal of the Indian Mathematical Society for more than a decade. He also served as the President of the prestigious Society. During his 34 year old association with
the TIFR, Ramanathan also visited other centres of mathematics, like the Institute for
Advanced Study at Princeton, the University of Missouri at St. Louis and the University of
Alberta at Edmonton, to name a few. He had authored two books.  He was also an influential mentor, guiding several distinguished students, many of whom continued to work in the Math Department of TIFR. Below listed are some of his doctoral students. The list has been obtained from an excellent review article (based on a lecture given at IMSc, Chennai) \cite{DP} by the eminent Indian mathematician, Prof. Dipendra Prasad.
\begin{itemize}
    \item  C . P Ramanujam - He was a talented mathematician, working in the areas of number theory and algebraic geometry. Like his namesake Ramanujan, he too died early after suffering from illness.
    \item  K. Ramachandra - He was a leading figure of Analytic, Algebraic and Transcendental Number Theory in India, well known for his work on the Riemann Zeta function and the six exponential theorem. His doctoral students include influential Number Theorists R. Balasubramanian and T.N Shorey.
    \item  S. Raghavan - Like his mentor KGR, he was also a winner of the prestigious Shanti Swarup Bhatnagar Award and was interested in the mathematics of Ramanujan. During his time at TIFR, he was the Dean of the School of Mathematics. From 1978-1982, he served as the Editor of the prestigious journal Proceedings(Mathematical Sciences), published by the Indian Academy of Sciences. It is also worth mentioning that one of his papers was the first from TIFR to appear in the Annals of Mathematics \cite{SR}. 
    \item Neela S. Rege - Born on 7th November 1941 in Baroda, Gujarat, she completed her Bachelors and Masters from the from the University of Bombay in 1963, and subsequently earned her PhD from TIFR in 1968. She received the G. Rangildas Mathematics prize in 1961 and the P. Thackersey award in 1963 from the University of Bombay. Dr. Rege is member of many learned Societies, like the American Mathematical Society, London Mathematical Society and Ramanujan Mathematical Society. During 1975-77, she was a Visiting Professor at U. dakar, Senegal. Since 1989, she is a Senior Lecturer at P.N.G University of Technology, Lae, Morobe, Papua New Guinea. \cite{NR}
    \item  Sunder Lal - Born on 11th October,1934 he received his MA (Mathematics) from Panjab University in 1958 and PhD degree from TIFR under the supervision of KGR in 1965, working in the area of one variable Modular forms. He joined the Department of Mathematics at Punjab University as a Reader on 13th July, 1966, and served there for about three decades, retiring in October 1994. \cite{SL}
    \item  V.C Nanda - Like Dr. Lal, he was also a Professor of Mathematics at Punjab University.
    \item S.S Rangachari - He was an eminent Number Theorist, working at the Tata Institute of Fundamental Research in Mumbai. He was interested in the works of Ramanujan. Post retirement, he had settled down in the United States of America.
\end{itemize}
\section{Later Life}
Apart from academics, Ramanathan was deeply interested in music and was a
singer himself. He was also fond of English, Tamil and Telugu literature. When he was at
Princeton, for about a couple of years, his neighbour was none other than Albert Einstein. He used to sing Carnatic songs to Einstein for his entertainment, especially songs of Tyagaraja. As a person, Ramanathan was a simple and
kind person, who could be easily approached by students. He was married to
Jayalakshmi Ramanathan, and the couple had 2 sons, Ananth and Mohan and 4 grandchildren Aparna, Kavitha, Anita and Nikil. In
December 1985, Ramanathan retired from the Tata Institute in Mumbai (erstwhile Bombay). During his last few years,
he suffered from Parkinson’s disease, and also underwent a cerebral surgery. After
suffering for a long time, Ramanathan passed away on May 10, 1992 at the age of 72,
marking an end to a very significant chapter in the history of Indian mathematics. The Proceedings (Mathematical Sciences) published a special issue in his honour, titled \enquote{K.G Ramanathan Memorial Issue}. Renowned mathematicians from all over the world (including some of KGR's own students and grandstudents), including Roger Heath Brown FRS, Raman Parimala, Don Zagier, R.P Bambah, R Balasubramanian, Bruce C Berndt, Richard Askey, George C Andrews contributed papers, and the issue received a review from another leading Indian mathematician, M.S Raghunathan FRS \cite{MSR}. Thus, K.G Ramanathan was a luminous figure in Indian mathematics post independence, and played a decisive role in bringing us back to the international mathematical map. A front ranked mathematician, a legendary institution builder and a great human being, KGR’s legacy continues even till
this date, and would definitely inspire the future generations of Indian students to take
up a life long career in mathematics.

\section{Acknowledgement}

The author would like to thank the anonymous referee for a careful reading of the manuscript and for many valuable comments. He is also grateful to the IACS Library (Cenenary Building) for providing him with a copy of the last reference.

\end{document}